\documentclass[a4paper,12pt,reqno]{amsart}
\usepackage{amssymb,amsthm,amsmath} 
\usepackage[english]{babel}
\usepackage{fullpage} 
\usepackage{hyperref}
\title[Noncommutative principal fibrations]
{Principal fibrations over noncommutative spheres}
\date{19 April 2018}

\author{Michel Dubois-Violette, Xiao Han, Giovanni Landi}
\address[]{\textit{Michel Dubois-Violette}
\newline \indent
Universit\'e Paris XI,
B\^atiment 210, F-91 405 Orsay Cedex}
\email{Michel.Dubois-Violette@u-psud.fr}
\address[]{\textit{Xiao Han}
\newline \indent SISSA,
via Bonomea 265, 34136 Trieste, Italy}
\email{xiao.han@sissa.it}
\address[]{\textit{Giovanni Landi}
\newline \indent
Universit\`a di Trieste,
Via A. Valerio, 12/1, 34127  Trieste, Italy
\newline \indent
and INFN, Trieste, Italy}
\email{landi@units.it}

\numberwithin{equation}{section}
\newtheorem{theo}{Theorem}[section]
\newtheorem{lemm}[theo]{Lemma}
\newtheorem{prop}[theo]{Proposition}

\newtheorem{rema}[theo]{Remark}

\newcommand{\ii}{{\mathrm{i} }}

\newcommand{\nn}{\nonumber}

\newcommand{\ca}{\mathcal{A}}

\newcommand{\ccR}{\mathcal{R}}
\newcommand{\IC}{{\mathbb C}}
\newcommand{\IH}{{\mathbb H}}
\newcommand{\II}{\mbox{\rm 1\hspace {-.6em} l}}
\newcommand{\IM}{{\mathbb M}}
\newcommand{\IN}{{\mathbb N}}
\newcommand{\IP}{{\mathbb P}}
\newcommand{\IR}{{\mathbb R}}

\newcommand{\IS}{{\mathbb S}}
\newcommand{\IT}{{\mathbb T}}
\newcommand{\IZ}{{\mathbb Z}}

\newcommand{\car}{{\ca_R}}

\DeclareMathOperator{\Mat}{Mat}

\DeclareMathOperator{\SO}{SO}
\DeclareMathOperator{\SU}{SU}
\DeclareMathOperator{\su}{su}
\DeclareMathOperator{\U}{U}
\DeclareMathOperator{\tr}{tr}
\DeclareMathOperator{\ccc}{ch}
\def\ket#1{\left | #1 \right\rangle}  %
\def\bra#1{\left\langle #1\right |}  %
\def\hs#1#2{\left\langle #1,#2\right\rangle}  %
\newcommand{\wt}{\widetilde}

\newcommand{\beqa}{\begin{align}}
\newcommand{\eeqa}{\end{align}}
\newcommand{\beq}{\begin{equation}}
\newcommand{\eeq}{\end{equation}}

\begin{document}


\begin{abstract}
We present examples of noncommutative four-spheres that are base spaces of $\SU(2)$-principal bundles with noncommutative seven-spheres as total spaces. The noncommutative coordinate algebras of the four-spheres are generated by the entries of a projection which is invariant under the action of $\SU(2)$. We give conditions for the components of the Connes--Chern character of the projection to vanish but the second (the top) one. The latter is then a non zero Hochschild cycle that plays the role of the volume form for the noncommutative four-spheres. 
\end{abstract}

\maketitle

\tableofcontents
\parskip = .9 ex

\thispagestyle{empty}

\section{Introduction and motivations}

The duality between spaces and algebras of functions on the spaces is at the basis of noncommutative geometry. One gives up the commutativity of the algebras of functions while replacing them by appropriate classes of noncommutative associative algebras which are considered as ``algebras of functions" on (virtual) ``noncommutative spaces".
 
For instance, one may consider noncommutative associative algebras generated by coordinates functions that satisfy relations other than the commutation between them, thus generalizing the polynomial algebras and defining thereby noncommutative vector spaces.

In this context, in the papers \cite{dvl17} and  \cite{dvl18} there were defined noncommutative finite-dimensional Euclidean spaces and noncommutative products of them. These ``spaces" were given in the general framework of the theory of regular algebras, which are a natural noncommutative generalizations of the algebras of polynomials. There were also defined noncommutative spheres  and noncommutative product of spheres. These are examples of noncommutative spherical manifolds related to the vanishing of suitable Connes--Chern classes of projections or unitaries in the sense of 
the work \cite{cl01} and \cite{cdv02}. 

In this paper we go one step further and consider actions of (classical) groups on noncommutative spheres and quotient thereof.  In particular we present examples of noncommutative four-spheres that are base spaces of $\SU(2)$-principal bundles with noncommutative seven-spheres as total spaces. This means that the noncommutative algebra of coordinate functions on the four-sphere is identified as a subalgebra of invariant elements for an action of the classical group $\SU(2)$ on the noncommutative algebra of coordinate functions on the seven-sphere. Conditions for these to qualify as noncommutative principal bundles are satisfied. The four-sphere algebra is generated by the entries of a projection which yields a noncommutative vector bundle over the sphere. Under suitable conditions the components of the Connes--Chern character of the projection vanish but the second (the top) one. The latter is then a non zero Hochschild cycle that plays the role of the volume form for the noncommutative four-sphere. 

The plan of the paper is organized as follows. In Section \ref{s:qa} we recall from \cite{dvl17} and \cite{dvl18} 
results on the quadratic algebras that we need and some of the solutions for noncommutative spheres 
$\IS^{7}_{R}$ that we use later on (here $R$ is a matrix of deformation parameters). In Section \ref{s:pb}, out of the functions on the seven-sphere $\IS^{7}_{R}$ we construct a projection in a matrix algebra over these functions, whose entries are invariant for a right action of $\SU(2)$ and thus generate a subalgebra that we identify as the coordinate algebra of a four-sphere $\IS^{4}_{R}$. We also show that this algebra inclusion is a noncommutative principal bundles with classical structure group $\SU(2)$. In Section \ref{s:star} 
the $*$-structure on the algebra of functions of $\IS^{4}_{R}$ is related to the vanishing of a component of the Connes--Chern character of the projection (and of a related unitary), a fact that put some restriction on the possible deformation matrices $R$, but makes the spheres examples of noncommutative spherical manifolds \cite{cl01}, \cite{cdv02}. We also exhibit explicit families of noncommutative 
four-sphere algebras. The Appendix~\ref{se:qpb} is devoted to some very basic notions of noncommutative principal bundles, while the Appendix~\ref{se:ccvf} relates Connes--Chern characters of idempotents and unitaries to Hochschild cycles and noncommutative volume forms.

\subsubsection*{Notation}
We use Einstein convention of summing over repeated up-down indices.
An algebra is always an associative algebra and a graded algebra is meant to be $\IN$-graded.


\section{A family of quadratic algebras}\label{s:qa}

\subsection{General definitions and properties}
In \cite{dvl17} and \cite{dvl18} there were
considered complex algebras $\car$ generated by two sets of hermitian elements
$x=(x_1, x_2)=(x_1^\lambda, x_2^\alpha)$, with $\lambda \in\{1, \dots, N_1\}$
and $\alpha \in\{1, \dots, N_2\}$, subject to relations
\begin{align}\label{alde}
& x_1^\lambda x_1^\mu = x_1^\mu x_1^\lambda \ ,
\qquad x_2^\alpha x_2^\beta = x_2^\beta x_2^\alpha , \nn \\
& x_1^\lambda x_2^\alpha = R^{\lambda\alpha}_{\beta\mu} \, x_2^\beta x_1^\mu \ , \qquad 
x^\alpha_2 x^\lambda_1= \overline{R}^{\lambda\alpha}_{\beta\mu} \> x^\mu_1 x^\beta_2
\end{align}
for a  `matrix' $(R^{\lambda\alpha}_{\beta\mu})$.
Here $\overline{R}^{\lambda\alpha}_{\beta\mu} \in \IC$ is the complex conjugates of the $R^{\lambda\alpha}_{\beta\mu}\in \IC$. The class of relevant matrices $R$ was defined by a series of conditions that we recall momentarily.

The quadratic complex algebra $\ca_R$ is a graded algebra $\ca_R=\oplus_{n\in\mathbb N}(\ca_R)_n$ which is connected, that is $(\ca_R)_0=\mathbb C\II$. Moreover, the quadratic relations \eqref{alde} of $\ca_R$ imply that there is a unique structure of $\ast$-algebra on $\ca_R$ for which the $x^\lambda_1$ $(\lambda\in \{1,\dots,N_1\})$ and the $x^\alpha_2$ $(\alpha\in \{1,\dots, N_2\})$ are hermitian, $x^\lambda_1=(x^\lambda_1)^\ast$ and $x^\alpha_2=(x^\alpha_2)^\ast$. This structure is graded in the sense that one has $f^\ast \in(\ca_R)_n\Leftrightarrow f\in (\ca_R)_n$ and $\ca_R$ is the quadratic $\ast$-algebra generated by the hermitian elements  $x^\lambda_1$ and $x^\alpha_2$ with the relations (\ref{alde}).

The $x^\lambda_1 x^a_1$ for $\lambda\leq \mu$ and the $x^\alpha_2 x^\beta_2$ for $\alpha\leq \beta$ are linearly independent in $(\ca_R)_2$ and generate $(\ca_R)_2$ together with the $x^\lambda_1 x^\alpha_2$. It is also natural to assume that the $x^\alpha_2 x^\lambda_1$ are independent which implies the equations
\begin{equation}
\overline{R}^{\lambda\alpha}_{\beta\mu} R^{\mu\beta}_{\gamma\nu}= \delta^\lambda_\nu \delta^\alpha_\gamma
\label{real}
\end{equation}
which in turn imply that the $x^\lambda_1 x^\alpha_2$ are also independent. Finally this implies in particular that the $x^\lambda_1 x^a_1$ with $\lambda\leq \mu$, the $x^\alpha_2 x^\beta_2$ with $\alpha\leq \beta$ and the $x^\nu_1 x^\gamma_2$ define a basis of $(\ca_R)_2$ while by definition the elements $x^\lambda_1$ and the $x^\alpha_2$ form a basis of $(\ca_R)_1$.


The classical (commutative) solution is given by
$$
(R_0)^{\lambda\alpha}_{\beta\mu}=\delta^\lambda_\mu  \delta^\alpha_\beta
$$
and $\ca_{R_0}$ is the coordinate algebra over the product $\IR^{N_1} \times \IR^{N_2}$.
Thus, the algebra $\car$ is though to define by duality the
noncommutative product of $\IR^{N_1} \times_R \IR^{N_2}$, that is $\car$ is the algebra
of coordinate functions on the noncommutative vector space $\IR^{N_1} \times_R \IR^{N_2}$.

If we collect together the coordinates, defining the
$x^a$ for $a\in \{1,2, \dots, N_1+N_2\}$ by $x^\lambda=x^\lambda_1$ and 
$x^{\alpha+N_1}=x^\alpha_2$, the relations \eqref{alde} with \eqref{real} can be written in the form
 \beq\label{aldeb}
x^a x^b = \ccR^{a \, b}_{c \, d} \,\, x^c x^d \ . 
\eeq
The $\ccR^{ab}_{cd}$ are the matrix elements of an endomorphism $\ccR$ of $(\ca_R)_1\otimes (\ca_R)_1$.
It follows from \eqref{real} that the $\ccR$ matrix  is involutive, that is
\beq\label{involr}
\ccR^2 = I \otimes I
\eeq
where $I$ is the identity mapping of $(\ca_R)_1$ onto itself. One 
next  imposes that the matrix $\ccR$ satisfies the Yang-Baxter equation
\beq\label{rvan}
(\ccR \otimes I) (I\otimes \ccR) (\ccR \otimes I) = (I \otimes \ccR)(\ccR \otimes I) 
(I \otimes \ccR) , 
\eeq
which then breaks in a series of conditions on the starting matrix $R^{\lambda\alpha}_{\beta\mu}$ in \eqref{alde}. 
Finally, additional conditions on the matrix $R^{\lambda\alpha}_{\beta\mu}$ comes by requiring that both quadratic elements
$(x_1)^2 = \sum^{N_1}_{\lambda=0}(x^\lambda_1)^2$ and
$(x_2)^2 = \sum^{N_2}_{\alpha=0}(x^\alpha_2)^2$ of $\car$ be central.

We refer to \cite{dvl17,dvl18} for a full analysis of the conditions that the matrix $R^{\lambda\alpha}_{\beta\mu}$ is required to satisfy. Here we  mention that the reality condition 
\eqref{real}, together with the requirement that the quadratic elements $(x_1)^2$ and 
$(x_2)^2$ be central,  leads to the symmetry conditions
\beq\label{eucl0}
R^{\lambda\beta}_{\alpha\mu}=R^{\mu\alpha}_{\beta\lambda}
=\overline{R}^{\mu\beta}_{\alpha\lambda}=(R^{-1})^{\beta\mu}_{\lambda\alpha} 
\eeq
as well as to the quadratic conditions
\beq\label{eucl1}
R^{\lambda\beta}_{\alpha\rho} R^{\rho\delta}_{\gamma\mu} = R^{\lambda\delta}_{\gamma\rho}R^{\rho \beta}_{\alpha\mu} \qquad \mbox{and} \qquad
R^{\lambda\beta}_{\gamma\nu} R^{\mu\gamma}_{\alpha\rho} =
R^{\mu\beta}_{\gamma\rho} R^{\lambda\gamma}_{\alpha\nu}
\eeq
which then (under \eqref{eucl0}) are equivalent to the cubic relations of the Yang--Baxter equations.
The general solution of these equations was given in \cite{dvl18} as follows.
By setting
$$
\widehat{R}^{\lambda\alpha}_{\mu\beta} = R^{\lambda\alpha}_{\beta\mu}
$$ 
for the endomorphism $\hat R = (\widehat{R}^{\lambda\alpha}_{\mu\beta})$ of $\mathbb R^{N_1}\otimes \mathbb R^{N_2}$ one has the representation
\beq
\label{ABCD}
\widehat{R}=\sum_r A_r\otimes B_r+ \ii \sum_a C_a\otimes D_a
\eeq
with the $A_r$ real symmetric $N_1 \times N_1$ matrices and the $B_r$ real symmetric $N_2 \times N_2$ matrices, both set taken to be linearly independent; 
and the $C_a$  real anti-symmetric $N_1 \times N_1$ matrices and the $D_a$ real anti-symmetric $N_2 \times N_2$ matrices (again both set taken to be linearly independent). Furthermore, they are such that 
\beq\label{AC}
[A_r,A_s]=0, \quad [A_r,C_a]=0, \quad [C_a,C_b]=0
\eeq
\beq\label{BD}
[B_r,B_s]=0, \quad [B_r,D_a]=0, \quad [D_a,D_b]=0
\eeq
for $r,s\in \{1,\dots, p\}$ and  $a,b\in\{1,\dots,q\}$,
with normalization condition 
\beq\label{NABCD}
\sum_{r,s} A_rA_s\otimes B_rB_s + \sum_{a,b} C_a C_b\otimes D_aD_b=\II_{N_1} \otimes \II_{N_2} ,
\eeq
a translation of the condition in \eqref{involr}

With the quadratic elements $(x_1)^2 = \sum^{N_1}_{\lambda=1}(x^\lambda_1)^2$ and 
$(x_2)^2 = \sum^{N_2}_{\alpha=1}(x^\alpha_2)^2$ of $\car$ being central, one may consider the quotient algebra
$$
\car / \left(  (x_1)^2-\II ,  (x_2)^2-\II \right)
$$
which defines by duality the noncommutative product 
$
\IS^{N_1-1} \times_R \IS^{N_2-1}
$
of the classical spheres $\IS^{N_1-1}$ and $\IS^{N_2-1}$. Indeed, for $R=R_0$, the above quotient is the restriction to $\IS^{N_1-1}\times \IS^{N_2-1}$ of the polynomial functions on 
$\IR^{N_1+N_2}$. 

Furthermore, with the central quadratic element $x^2  = \sum^{N_1+N_2}_{a=1} (x^a)^2 
= (x_1)^2 +  (x_2)^2$,
one may also consider the quotient of $\car$
$$
\car/ \left(  x^2 -\II \right).
$$
This defines (by duality) the noncommutative $(N_1+N_2-1)$-sphere 
$
\IS^{N_1+N_2-1}_R
$
shown in \cite{dvl18} to be a noncommutative spherical manifold in the sense of 
\cite{cl01} and \cite{cdv02}.

\subsection{Some quaternionic geometry}
When $N_1=N_2=4$, explicit solutions for the matrix $R^{\lambda \alpha}_{\beta \mu}$ 
were given in \cite{dvl17} and \cite{dvl18}
by using results on the geometry of quaternions.

The space of quaternions $\IH$ is identified with $\IR^4$ in the usual way:
\beq\label{qrid}
\IH \ni q = x^0 e_0 + x^1 e_1 + x^2 e_2 + x^3 e_3 \quad \longmapsto \quad x = (x^\mu) =
(x^0, x^1, x^2, x^3) \in \IR^4 .
\eeq
Here $e_0=1$ and the imaginary units $e_a$ obey the multiplication rule of the algebra $\IH$:
$$
e_a e_b = - \delta_{ab} + \sum_{c=1}^3\varepsilon_{abc} e_c .
$$
From this it follows an identification of the unit quaternions $\U_1(\IH) = \{q \in \IH \, | \, q \bar q = 1 \} $ with the euclidean three-sphere $\IS^3 = \{x \in \IR^4 \, ; \, || x ||^2 = \sum_\mu (x^\mu)^2 = 1\}$.

With the identification \eqref{qrid}, left and right multiplication of quaternions are represented
by matrices acting on $\IR^4$:
$$
L_{q'} q := q'q \quad \to \quad E^+_{q'} (x) \qquad \mbox{and} \qquad
R_{q'} q := qq' \quad \to \quad E^-_{q'} (x) .
$$
For $q$ a unit quaternion, both $E^+_{q}$ and $E^-_{q}$ are orthogonal matrices. In fact
the unit quaternions form a subgroup of the multiplicative group $\IH^*$ of non vanishing quaternions.
When restricting to these, one has then the identification
$$
\U_1(\IH) \simeq \SU(2),
$$
that is $E^+_{q}$ and $E^-_{q}$, for $q\in \U_1(\IH)$, are commuting $\SU(2)$ actions (each in the `defining representation')
on $\IR^4$, or together an action of $\SU(2)_L \times \SU(2)_R$ on $\IR^4$, with $L/R$ denoting left/right action.
This action is the adjoint one, an action of $\SO(4) = \SU(2)_L \times \SU(2)_R / \IZ_2$.

Let us denote $E^{\pm}_{a} = E^{\pm}_{e_a}$ for the imaginary units.
By definition one has that
\begin{align*}
E^{+}_{a} E^{-}_{b} = E^{-}_{b} E^{+}_{a} , \qquad  E^{\pm}_{a} E^{\pm}_{b} = - \delta_{ab} \II
\pm \sum_{c=1}^3\varepsilon_{abc} E^{\pm}_{c} .
\end{align*}
In the following, it will turn out to be more convenient to change a sign to the `right' matrices: we shall rather use matrices
$J^{+}_{a}:=E^{+}_{a}$ and $J^{-}_{a}:=-E^{-}_{a}$. For these one has
\begin{align*}
J^{+}_{a} J^{-}_{b} = J^{-}_{b} J^{+}_{a} , \qquad
J^{\pm}_{a} J^{\pm}_{b} = - \delta_{ab} \II + \sum_{c=1}^3\varepsilon_{abc} J^{\pm}_{c} ,
\end{align*}
that is the matrices $J^\pm_a$ are two copies of the quaternionic imaginary units.
Indeed for these $4 \times 4$ real matrices $J^\pm_a$ one can explicitly compute
\beq\label{jmatrix}
(J^{\pm}_{a})_{\mu\nu} = \mp ( \delta_{0\mu}  \delta_{a\nu} - \delta_{a\mu}  \delta_{0\nu})
+ \sum_{b,c=1}^3\varepsilon_{abc} \delta_{b\mu} \delta_{c\nu} .
\eeq

With the identification $\U_1(\IH) \simeq \SU(2)$, when acting on $\IR^4$, the matrices $J^\pm_a$ are a  representation of
the Lie algebra $\su(2)$ of $\SU(2)$,
or taken together a representation of $\su(2)_L \oplus \su(2)_R$.
For the standard positive definite scalar product on $\IR^4$, the six matrices $J^{\pm}_{a}$ are readily checked to be antisymmetric, ${}^t J^{\pm}_{a} = - J^{\pm}_{a}$, and one finds in addition that 
$$
-\tfrac{1}{4} \tr (J^{\pm}_{a} J^{\pm}_{b}) =\delta_{ab} .
$$
Then, the matrices $(J^{\pm}_{1}, J^{\pm}_{2}, J^{\pm}_{3})$ are canonically an orthonormal basis of 
$\Lambda^2_\pm \IR^4{}^* \simeq \IR^3$ considered as an oriented three-dimensional euclidean space with the orientation of this basis; mapping $J^{\pm}_{a} \to J^{\mp}_{a}$ amounts to exchange the orientation. 
On the other hand, the nine matrices $J^{+}_{a} J^{-}_{b}$ are an orthonormal basis for the space of symmetric 
trace-less $4 \times 4$ matrices.

\subsection{Noncommutative quaternionic tori and spheres}\label{se:q7s}
Referring to the above, we have explicit solutions for the deformation matrix in \eqref{ABCD}.
Firstly, with any vector ${\bf u} = (u^1, u^2, u^3) \in \IR^3$ we get antisymmetric matrices
$$ 
J^{+}_{{\bf u}} := u^1 J^{+}_{1} + u^2 J^{+}_{2} + u^3 J^{+}_{3} \qquad \mbox{or} \qquad
J^{-}_{{\bf u}} := u^1 J^{-}_{1} + u^2 J^{-}_{2} + u^3 J^{-}_{3} .
$$
With this notation, consider the matrix 
\beq\label{r8dquat0} 
R^{\lambda \alpha}_{\beta \mu} = u^0 \, \delta^\lambda_\mu \delta^\alpha_\beta +
\ii \, (J^+_{{\bf v}})^\lambda_\mu \, (J^+_{{\bf u}})^\alpha_\beta .
\eeq
Clearly, all the commutation relation \eqref{AC} and \eqref{BD} are satisfied. Thus, with this matrix, 
we define $\ca_R$ as the $\ast$-algebra generated by the hermitian elements $x^\lambda_1$ and $x^\alpha_2$, $\lambda,\alpha\in\{0,1,2,3\}$,  
with relations
\beq\label{ald4}
x^\lambda_1 x^\mu_1=x^\mu_1 x^\lambda_1,\qquad  x^\alpha_2x^\beta_2=x^\beta_2x^\alpha_2,\qquad  x^\lambda_1 x^\alpha_2=R^{\lambda\alpha}_{\beta\mu} x^\beta_2 x^\mu_1
\eeq
But using the action of $\SO(3)$ one can always rotate ${\bf v}$ to a fixed direction
$\widehat{{\bf u}}$,
and in this case the resulting matrix $R$ has parameters $u^0\in\IR$ and ${\bf u}\in\IR^3$
constrained to
$$
(u^0)^2 + {\bf u}^2 = 1,
$$
that is they make up a three-dimensional sphere $\IS^3$.  There is in fact a residual `gauge' freedom in that one can use a rotation around the direction $\widehat{{\bf u}}$ to remove one component of the vector ${\bf u}$. Thus if $\widehat{{\bf u}}_1$ and $\widehat{{\bf u}}_2$ are two orthogonal unit vectors (say in the canonical basis), we get families of noncommutative spaces determined by the matrices
\beq\label{r8dquat}
R^{\lambda \alpha}_{\beta \mu} = u^0 \delta^\lambda_\mu \delta^\alpha_\beta +
\ii \, (J^+_{1})^\lambda_\mu \, (u^1 J^+_{1} + u^2 J^+_{2})^\alpha_\beta ,
\eeq
and parameters constrained by a two-dimensional sphere $\IP^1(\IC) = \IS^3 / \IS^1 = \IS^2$ being
$$
(u^0)^2 + (u^1)^2 + (u^2)^2 = 1 .
$$

These constructions lead to natural quaternionic generalisations of the toric
four-dimensional noncommutative spaces described in \cite{cl01} for which the space of deformation parameter is $\IP^1(\IR) = \IS^1 / \IZ_2 = \IS^1$.

Indeed, in parallel to the complex case were there is an action of the classical torus $\IT^2$, 
there is now an action of the classical quaternionic torus  $T^2_\IH=\U_1(\IH)\times \U_1(\IH)=\IS^3\times \IS^3=\SU(2)\times \SU(2)$ by $\ast$-automorphisms of the  algebra $\ca_R$ given as follows. 
%

%
%

In view of the commutations of the $J^-_a$ with the $J^+_b$ for $a,b\in \{1,2,3\}$, the mappings $x_1\mapsto J^-_a x_1$, $x_2\mapsto J^-_bx_2$ for $a,b\in \{1,2,3\}$ leave the relations \eqref{ald4} of $\ca_R$ invariant and thus define $\ast$-automorphisms of the $\ast$-algebra $\ca_R$.  
By setting $q=q^0+q^a e_a\in \IH^\ast\mapsto q^0\II +q^aJ^-_a$ with obvious conventions, one has from last section (right quaternionic) actions $x_1\mapsto (q^0_1\II +q^a_1 J^-_a)x_1$ and  $x_2\mapsto (q^0_2\II +q^a_2J^-_a) x_2$
of the multiplicative group $\IH^\ast\times \IH^\ast$ as automorphisms of the $\ast$-algebra $\ca_R$, 

This induces an action of $\U_1(\IH)\times \U_1(\IH)$ on $\ca_R$ by restriction to the $q\in \U_1(\IH)$ which passes to the quotient by the ideal generated by the two central elements $(x_1)^2 = \sum_\lambda (x^\lambda_1)^2$, $(x_2)^2=\sum_\alpha (x^\alpha_2)^2$ and defines an action of the classical quaternionic torus  $\U_1(\IH)\times \U_1(\IH)$ by $\ast$-automorphisms of the coordinate algebra 
$$
\ca((T^2_{\IH})_R) = \car / \left(  (x_1)^2-\II ,  (x_2)^2-\II \right)
$$
of the ``noncommutative" quaternionic torus $(T^2_{\IH})_R$.    The action also passes to 
the quotient by the ideal generated by the central element $(x_1)^2 + (x_2)^2$ and defines an action of the classical quaternionic torus  $\U_1(\IH)\times \U_1(\IH)$ by $\ast$-automorphisms of the coordinate algebra 
$$
\ca(\IS^{7}_{R}) = \car / \left(  (x_1)^2 + (x_2)^2 -\II \right)
$$
of a noncommutative seven-sphere $\IS^{7}_{R}$. As we shall see in what follows, when restricting to the diagonal action of $\U_1(\IH) \subset \U_1(\IH)\times \U_1(\IH)$ on $\ca(\IS^{7}_{R})$
will result into a $\SU(2)$-principal bundles $\IS^{7}_{R} \to \IS^{4}_{R}$ on a noncommutative four-sphere.

\section{Principal fibrations}\label{s:pb}
We are going to define natural  $\SU(2)$-principal bundles $\IS^{7}_{R} \to \IS^{4}_{R}$ in the `dual' sense of a coordinate algebra $\ca(\IS^{4}_{R})$ on a four-sphere that is identified 
as the invariant subalgebra of the coordinate algebra $\ca(\IS^{7}_{R})$ on a seven-sphere, for an action of the group $\SU(2)$.

\subsection{A canonical projection}
In parallel with \eqref{qrid} consider the two quaternions
$$
x_1 = x_1^\mu e_\mu , \quad x_2 = x_2^\alpha e_\alpha ,
$$
with commutation relations among their components governed by a matrix $R^{\lambda\alpha}_{\beta\mu}$ as in \eqref{alde}.
Then, when restricting to the sphere $\IS^{7}_{R}$ the vector-valued function
\beq\label{ket}
\ket{\psi} = \begin{pmatrix} x_2 \\ x_1 \end{pmatrix}
\eeq
has norm
$$
\hs{\psi}{\psi} = ||x_1||^2 + ||x_2||^2 = \II
$$
and thus we get a projection
\beq\label{pro}
p = \ket{\psi}\bra{\psi} = \begin{pmatrix} x_2 x^*_2  & x_2 x^*_1 \\
x_1 x^*_2 & x_1 x^*_1
\end{pmatrix} ,
\eeq
that is $p=p^*=p^2$. Define coordinate functions $Y= Y^0 e_0 + Y^k e_k$ and $Y^4$ by
\begin{align}\label{ydef}
Y^4 &= ||x_2||^2 - ||x_1||^2  \qquad \mbox{and} \qquad
\tfrac{1}{2} Y = x_2 x^*_1
\end{align}
so that the projection \eqref{pro} is written as
\beq\label{pro1}
p = \ket{\psi}\bra{\psi} = \tfrac{1}{2} \begin{pmatrix} 1 + Y^4 & Y \\
Y^* & 1 - Y^4
\end{pmatrix} .
\eeq
The condition $p^2=p$ leads to
\begin{align}\label{cond0}
Y Y^* + (Y^4)^2 = 1 \quad &\mbox{and} \quad Y^* Y + (Y^4)^2 = 1 \\
Y Y^4 = Y^4 Y \quad &\mbox{and} \quad Y^* Y^4 = Y^4 Y^*  \label{cond00}.
\end{align}
Thus the coordinate function $Y^4$ is central while comparing the first two conditions requires
$Y Y^*=Y^* Y$ and that this is a (central) multiple of the identity. A direct computation translates these to the conditions
\begin{align}
-(Y^{0*} Y^k - Y^{k*} Y^0 ) + \varepsilon_{kmn} Y^{m*} Y^{n} &= 0  \label{4sp1},  \\
Y^0 Y^{k*} - Y^k Y^{0*} + \varepsilon_{kmn} Y^{m} Y^{n*} & = 0  \label{4sp2}
\end{align}
for $k,r,m=1,2,3$ and totally antisymmetric tensor $\varepsilon_{krm}$, together with 
\beq\label{4sp}
\sum_{\mu=0}^{3} (Y^{\mu*} Y^{\mu} - Y^{\mu}Y^{\mu*} ) = 0 .
\eeq
Then condition \eqref{cond0} reduces to a four-sphere relation
\beq\label{4sr}
\sum_{\mu=0}^{3} Y^{\mu*} Y^{\mu} + (Y^4)^2 = 1 = \sum_{\mu=0}^{3} Y^{\mu} Y^{\mu*} + (Y^4)^2 .
\eeq
Being $Y^4$ central, these relations also give that both $\sum_{\mu=0}^{3} Y^{\mu*} Y^{\mu}$ and $\sum_{\mu=0}^{3} Y^{\mu} Y^{\mu*}$ are central as well.
In view of the relations \eqref{4sr}, the elements $Y^{\mu}$ generate the $*$-algebra
$\ca(\IS^{4}_R)$ of a four-sphere $\IS^{4}_R$. This four-sphere $\IS^{4}_R$ is the suspension 
(by the central element $Y^4$), of a three-sphere $\IS^{3}_R$ obtained by reducing \eqref{4sr} to 
\beq\label{3sr}
\sum_{\mu=0}^{3} Y^{\mu*} Y^{\mu} = 1 = \sum_{\mu=0}^{3} Y^{\mu} Y^{\mu*} .
\eeq

\begin{rema}
Up to the change $Y^0 \mapsto -Y^0$, the relations \eqref{4sp1}-\eqref{4sp2} are the same as the relations (2.4)-(2.6) of \cite{cdv02}.  
\end{rema}

Clearly, the coordinate function $Y^4$ is hermitian. On the other hand, as we shall see, the coordinate functions $Y^{\mu*}$, $\mu=0,1,2,3$, while not hermitian, are not independent from the $Y^{\mu}$'s with
the explicit dependence determined by the matrix $R^{\lambda\alpha}_{\beta\mu}$.
In fact the commutation relations \eqref{4sp1}-\eqref{4sp2} are not additional relations but they are determined by the $R^{\lambda\alpha}_{\beta\mu}$ which gives the commutation relations among the starting $x$'s.

\begin{lemm}
With the definitions in \eqref{ydef} it holds that
\beq\label{y0}
\tfrac{1}{2} Y^0 = \sum_{\mu=0}^{3} x_2^{\mu} x_1^{\mu} , \qquad
\tfrac{1}{2} Y^k = x_2^{k} x_1^{0} - x_2^{0} x_1^{k} - \varepsilon_{knm} x_2^{n} x_1^{m}
\eeq
and
\beq\label{y*0}
\tfrac{1}{2} Y^{0*} = \sum_{\mu=0}^{3} x_1^{\mu} x_2^{\mu} , \qquad
\tfrac{1}{2} Y^{k*} = x_1^{0} x_2^{k} - x_1^{k} x_2^{0} + \varepsilon_{knm} x_1^{n} x_2^{m} \, .
\eeq
\begin{proof}
A direct computation.
\end{proof}
\end{lemm}

\subsection{Noncommutative $\SU(2)$-principal bundles}
As mentioned, due to the relations \eqref{4sr}, the elements $Y^{\mu}$ generate the $*$-algebra
$\ca(\IS^{4}_R)$ of a four-sphere $\IS^{4}_R$. The algebra inclusion $\ca(\IS^{4}_R) \hookrightarrow  \ca(\IS^{7}_R)$ is a principal $\SU(2)$ bundle in the following sense.

With $\ket{\psi}$ the vector-valued function in \eqref{ket}, let the action of a unit quaternion
$w \in \U_1(\IH) \simeq \SU(2)$ on $\IS^{7}_R$ be obtained from the following action on the generators:
\beq\label{actsu0}
\alpha_w(\ket{\psi}) = \ket{\psi}w = \begin{pmatrix} x_2 w \\ x_1 w \end{pmatrix}.
\eeq
Clearly, the projection $p$ and then the algebra $\ca(\IS^{4}_R)$ are invariant for this action. 

On the other hand, in general the action \eqref{actsu0} does not preserve the commutation relations of the $\IS^{7}_R$ and thus not results 
into an action by $\ast$-automorphisms of the coordinate algebra 
$\ca(\IS^{7}_{R})$. 
Let us assume this is the case, that is the action preserves the commutation relations and postpone to later on the study of deformations that meet this condition.
If we denote $H=\ca(\SU(2))$ we have dually a co-action $\delta$ of $H$ by $\ast$-automorphisms on $\ca(\IS^{7}_R)$ with algebra of co-invariant element again the subalgebra $\ca(\IS^{4}_R)$.
\begin{prop} 
Let the action of $\SU(2)$ in \eqref{actsu0} be 
by $\ast$-automorphisms of the coordinate algebra $\ca(\IS^{7}_{R})$
and let $H=\ca(\SU(2))$. Then the canonical map
$$
\chi :  \ca(\IS^{7}_R) \otimes_{\ca(\IS^{4}_R)} \ca(\IS^{7}_R) \to \ca(\IS^{7}_R) \otimes H  ,
\qquad \chi(p' \otimes p) = p' \delta(p)
$$
is bijective.
\begin{proof}
With $\ket{\psi}$ as in \eqref{ket}, one has that
$$
\chi \big( \bra{\psi} \otimes_{\ca(\IS^{4}_R)} \ket{\psi} \big) =
\bra{\psi} \delta (\ket{\psi}) =  \hs{\psi}{\psi} \otimes w = \II \otimes w .
$$
As described in the Appendix~\ref{se:qpb} this is enough for the surjectivity of $\chi$. Moreover, being the structure group classical (and semisimple), again as mentioned in that appendix, surjectivity is equivalent to bijectivity (see \cite{lvs05} for a similar construction and example).
\end{proof}
\end{prop}

\subsection{Volume forms}\label{s:vf}

We have already observed that the unit radius conditions in \eqref{cond0} requires
that the `quaternion' $Y= Y^0 e_0 + Y^k e_k$ be such that 
$Y Y^*=Y^* Y \in \II_2 \otimes \ca(\IS^{4}_R)$ (in fact be in the centre of $\ca(\IS^{4}_R)$.  
An important role is played by the components of the Connes--Chern character in cyclic homology of $Y$ 
(see Appendix~\ref{se:ccvf}),
\beq\label{vacc1}
\ccc_{\frac{1}{2}}(Y) = \left< Y \otimes Y^* - Y^* \otimes Y \right> 
\eeq 
and  
\beq\label{vacc2}
\ccc_{\frac{3}{2}}(Y) = \left< Y \otimes Y^* \otimes Y \otimes Y^* -  Y^* \otimes Y \otimes Y^* \otimes Y \right> .
\eeq
Here $\left<  ~\cdot~ \right>$ indicates the partial matrix trace over $\IM_2({\IC})$, thinking of  $\IH$ as a subset of 
$\IM_2({\IC})$. 
From the appendix we know that when $\ccc_{\frac{1}{2}}(Y) = 0$ the element  $\ccc_{\frac{3}{2}}(Y)$ is a Hochschild cycle which gives a volume element for the three-sphere $\IS^{3}_R$ obtained by the unitarity 
conditions $Y Y^*=Y^* Y = \II_2$. 

On the other hand, for the projection $p$ in \eqref{pro1} one has at once that
$$
\ccc_0 (p) = \left< ( p - \tfrac{1}{2} \, \II) \right> = 0. 
$$
Moreover, being the four-sphere $\IS^{4}_R$ the suspension by the central element $Y^4$ of the three-sphere $\IS^{3}_R$, 
the vanishing $\ccc_{\frac{1}{2}}(Y) = 0$ would also imply the vanishing of the component $\ccc_1 (p)$ (cf. \cite[Theorem 2]{cdv02}), where
$$
\ccc_1 (p) = \left< ( p - \tfrac{1}{2} \, \II) \otimes p^{\otimes 2} \right>.
$$
Then, similarly to before, the element  $\ccc_2 (p) = \left< ( p - \tfrac{1}{2} \, \II) \otimes p^{\otimes 4} \right>$ 
is a Hochschild cycle which gives a volume element for the four-sphere $\IS^{4}_R$.

We see that the vanishing $\ccc_{\frac{1}{2}}(Y) = 0$ makes both the sphere $\IS^{3}_R$ 
as well as its suspension sphere $\IS^{4}_R$, noncommutative spherical manifolds in the sense of \cite{cl01}, \cite{cdv02}. 

\subsection{An analysis of the $*$-structure}\label{s:star}

Due to relation \eqref{4sp}, one expect the elements $Y^{a*}$, $a=0,1,2,3$, to be expressed in terms of the elements $Y^a$, $a=0,1,2,3$. This fact requires conditions on the possible deformation matrix $R$, while giving nicer properties for both spheres $\IS^{4}_R$ and $\IS^{3}_R$. Indeed, they becomes spherical manifolds as mentioned at the end of previous section.
A direct computation shows that the vanishing $\ccc_{\frac{1}{2}}(Y) = 0$ 
is equivalent to the condition
\beq\label{vanzz}
\sum_{\mu=0}^{3} (Y^{\mu*} \otimes Y^{\mu} - Y^{\mu} \otimes Y^{\mu*} ) = 0 .
\eeq
One has then the following 
\begin{lemm}\cite[Lemma 2]{cdv02}
The condition \eqref{vanzz} is satisfied, that is $\ccc_{\frac{1}{2}}(Y)=0$, 
if and only if there is a symmetric unitary matrix $\Lambda\in\IM_4(\IC)$ such that
\beq\label{lambda0}
Y^{\mu*} = \Lambda^\mu{}_\nu Y^\nu \qquad \mu,\nu=0,1,2,3 \, .
\eeq
\end{lemm}
%
%
\noindent
In turn, the condition $\ccc_{\frac{1}{2}}(Y) = 0$ is left unchanged by a linear change in generators as
\beq
Y^\mu \mapsto u \, S^\mu{}_\nu \, Y^\nu
\eeq
with $u\in\U(1)$ and $S\in\SO(4)$ a real rotation. Under this transformation,
the symmetric unitary matrix $\Lambda$ in \eqref{lambda0} transforms as
\beq
\Lambda \mapsto u^2 \, S^t \, \Lambda S .
\eeq
Then it can be diagonalized by a real rotation $S$ and with a further normalization
(by a factor $u\in\U(1)$) it can alway be put in the form
\beq\label{lambda}
\Lambda = \begin{pmatrix}
1 & 0                    & 0 &0 \\
0 & e^{\ii \theta_1} & 0 &0 \\
0 & 0 & e^{\ii \theta_2} & 0 \\
0 & 0 & 0 & e^{\ii \theta_3}
\end{pmatrix} ,
\eeq
for suitable angles $\theta_1,\theta_2,\theta_3$ (see \cite[\S2]{cdv02}).

%

\section{The quaternionic family of four-spheres}
Let us now consider the quaternionic deformations mentioned in section \ref{se:q7s} governed by the deformation matrix in \eqref{r8dquat} and in particular the noncommutative seven-sphere.  
As we have seen, there is a compatible action of $\U_1(\IH)\times \U_1(\IH)$ by $\ast$-automorphisms of the corresponding coordinate algebra.
Then for the action in \eqref{actsu0} we may take the diagonal action by $w=w^0+w^a e_a \in \U_1(\IH) \mapsto w^0\II +w^aJ^-_a\in\SU(2)$, written 
explicitly on generators as  $x_1\mapsto (w^0\II +w^aJ^-_a)x_1$ and  $x_2\mapsto (w^0\II +w^aJ^-_a) x_2$. 

\begin{prop}
Given the commutation relations for the $x$'s for the matrix \eqref{r8dquat}, one 
has $Y^{\mu*} = \Lambda^\mu{}_\nu Y^\nu$ for
$\Lambda\in\IM_4(\IC)$ a symmetric unitary matrix given explicitly by:
\beq\label{qua*1}
\begin{pmatrix}
Y^{0*} \\ Y^{3*}
\end{pmatrix}
=
\begin{pmatrix}
u^0 + \ii u^1 &  \ii u^2 \\
\ii u^2 & u^0 - \ii u^1
\end{pmatrix}
\begin{pmatrix}
Y^0 \\ Y^3
\end{pmatrix}
\eeq
and
\beq\label{tor*2}
\begin{pmatrix}
Y^{1*} \\ Y^{2*}
\end{pmatrix}
=
\begin{pmatrix}
u^0 + \ii u^1 &  \ii u^2 \\
\ii u^2 & u^0 - \ii u^1
\end{pmatrix}
\begin{pmatrix}
Y^1 \\ Y^2
\end{pmatrix} .
\eeq

\begin{proof}
A comparison with the matrices \eqref{jmatrix} shows that
$\tfrac{1}{2} Y^a = (J_+^a)_{\alpha\lambda} \, x_2^{\alpha} x_1^{\lambda}$
for $a=0,1,2,3$, with $J_+^0=\II$. This allows one to write
$\tfrac{1}{2} Y^{a*} = (J_+^a)_{\alpha\lambda} \, x_1^{\lambda} x_2^{\alpha}
= (J_+^a)_{\alpha\lambda} \, R^{\lambda \alpha}_{\beta \mu} \, x_2^\beta x_1^\mu$.
Then, 
with the $R$ matrix in \eqref{r8dquat}, a direct computation of \eqref{y*0} yields
\begin{align*}
\tfrac{1}{2} Y^{a*} &= (J^+_a)_{\alpha\lambda} \, R^{\lambda \alpha}_{\beta \mu} \, x_2^\beta x_1^\mu \\
&= \left( u^0 J^+_{a} -\ii \, u^1 J^+_{1} J^+_{a} J^+_{1} -\ii \, u^2 J^+_{2} J^+_{a} J^+_{1}
\right)_{\beta \mu}  x_2^\beta x_1^\mu \\
& = \tfrac{1}{2} \left( u^0 Y^a -\ii \, u^1 J^+_{1} Y^a J^+_{1} -\ii \, u^2 J^+_{2} Y^a J^+_{1} \right),
\end{align*}
from which one gets the explicit expressions in \eqref{qua*1} and \eqref{tor*2}.
\end{proof}
\end{prop}
With the $*$-structure of the previous proposition one shows that none of the generators is normal, that is 
$Y^{\mu*} Y^{\mu} \not= Y^{\mu} Y^{\mu*}$ while the condition \eqref{4sp} is automatically satisfied. On the other hand, 
the commutation relations \eqref{4sp1} and \eqref{4sp2} can be written as:
\begin{align*}
(u^0 + \ii u^1) (Y^1 Y^0 - Y^0 Y^1) + \ii u^2 (Y^1 Y^3 - Y^0 Y^2) & = 0 \\
(u^0 - \ii u^1) (Y^3 Y^2 - Y^2 Y^3)  + \ii u^2 (Y^3 Y^1 - Y^2 Y^0)  & = 0 \\
u^0 (Y^2 Y^0 - Y^0 Y^2)  - \ii u^1 (Y^1 Y^3 + Y^3 Y^1) + \ii u^2 ( Y^1 Y^0 - Y^3 Y^2 ) & = 0 \\
u^0 (Y^3 Y^1 - Y^1 Y^3)  - \ii u^1 (Y^0 Y^2 + Y^2 Y^0) + \ii u^2 ( Y^0 Y^1 - Y^2 Y^3 )  & = 0 \\
u^0 (Y^3 Y^0 - Y^0 Y^3)  + \ii u^1 (Y^1 Y^2 + Y^2 Y^1) + \ii u^2 \big( (Y^2)^2 - (Y^1)^2\big)  & = 0 \\
u^0 (Y^2 Y^1 - Y^1 Y^2)  + \ii u^1 (Y^0 Y^3 + Y^3 Y^0) + \ii u^2 \big( (Y^3)^2 - (Y^0)^2\big)  & = 0 .
\end{align*}

For the structure $\Lambda$ in \eqref{qua*1} and \eqref{tor*2}, the matrix
\beq
\Lambda' = \begin{pmatrix}
u^0 + \ii u^1 &  \ii u^2 \\
\ii u^2 & u^0 - \ii u^1
\end{pmatrix}
\eeq
being symmetric and unitary, can be diagonalized by a real rotation $S$: one finds eigenvalues 
$\lambda_{\pm} = u^0 \pm \ii \sqrt{(u^1)^2 + (u^2)^2} = u^0 \pm \ii \sqrt{1-(u^0)^2}$.  
With a further normalization by the factor 
$u^0 - \ii \sqrt{(u^1)^2 + (u^2)^2}  \in \U(1)$, the matrix $\Lambda'$ can be put in the form
\beq
\begin{pmatrix}
1 & 0 \\
0 & e^{\ii \theta}
\end{pmatrix} ,
\eeq
and a direct computation gives:
\beq
e^{\ii \theta} = \frac{ u^0 + \ii \sqrt{(u^1)^2 + (u^2)^2} }{ u^0 - \ii \sqrt{(u^1)^2 + (u^2)^2} } 
= \big( u^0 + \ii \sqrt{(u^1)^2 + (u^2)^2} \, \big)^2 .
\eeq
The sphere $\IS^{4}_R=\IS^{4}_\theta$ is then (isomorphic to) a $\theta$-deformation, as the one introduced in \cite{cl01}

\appendix

\section{Noncommutative principal bundles}\label{se:qpb}

Quantum (noncommutative) principal bundles we are interested in were introduced in \cite{bm93}. 
As a total space one considers an algebra $P$ and as structure group  
a Hopf algebra $H$. The algebra $P$ is a right $H$-comodule algebra
with right coaction $\Delta_{R}\,:\,P \mapsto P \otimes H$ 
(for which one uses Sweedler-like notation: $\Delta_R(p)= p_{(0)} \otimes p_{(1)}$).
The subalgebra of the right coinvariant elements, 
$$
B=P^H \,:=\,\{p\in P\,:\,\Delta_{R} p = p\otimes 1\} ,
$$
is the base space algebra of the bundle.
At the algebraic level the principality of the bundle is the requirement that the ``canonical map"
\beq\label{can}
\chi : P \otimes_B P \to P \otimes H; \quad p' \otimes_B p \mapsto p' \Delta_R(p) = p' p_{(0)} \otimes p_{(1)}
\eeq
is bijective. This is indeed the definition that the inclusion $B \hookrightarrow P$ be a Hopf-Galois extension \cite{sc90}. The canonical map is left $P$-linear and right $H$-colinear and is a morphism (an isomorphism  for Hopf-Galois extensions) of left $P$-modules and right $H$-comodules. It is also clear that $P$ is 
both a left and a right $B$-module.

For structure Hopf algebras $H$ which are cosemisimple and have bijective antipodes, Theorem~I of \cite{sc90} grants additional nice properties. In particular, the surjectivity of the canonical map implies its bijectivity.
Moreover, in order to prove surjectivity of $\chi$, it is enough to prove that for any generator $h$ of $H$, the element $1\otimes h$ is in the image of the canonical map. Indeed, if 
$\chi(g_k \otimes_{\mathcal{B}} g'_k) =1 \otimes g$ and $\chi(h_l \otimes_{\mathcal{B}} h'_l) =1 \otimes h$ for $g,h\in H$, 
then $\chi(g_k h_l \otimes_{\mathcal{B}} h_l' g_k')=g_k h_l \chi(1 \otimes_{\mathcal{B}} h_l' g_k')= 1\otimes h g$, using the fact that the canonical map restricted to $1\otimes_B P$ is a homomorphism. Extension to all of $P\otimes_{\mathcal{B}} P$ then follows from left $P$-linearity of $\chi$. It is also easy to write down an explicit expression for the inverse of the canonical map. Indeed, one has $\chi^{-1} (1 \otimes hg) = g_k h_l \otimes_B h_l' g_k'$ in the above notation so that the general form of the inverse follows again from left $P$-linearity. 

\section{Connes--Chern characters and volume forms}\label{se:ccvf}
Let $\ca$ be a unital algebra over $\IC$ and let 
$\wt{\ca} = \ca /  \IC \II $ be the quotient of $\ca$ by the scalar multiples of the unit $\II$.
Given an idempotent, 
$$
e = (e^{i}_{j}) \in \Mat_r (\ca) \, \qquad e^2 = e ,
$$
the component $\ccc_k (e)$ of the (reduced) Chern character of $e$ is the element of
$\ca \otimes (\wt{ \ca} )^{\otimes 2k}$,
given explicitly by the formula 
\begin{align}\label{B1}
\ccc_k (e) & = \lambda_k \left< ( e - \tfrac{1}{2} \, \II) \otimes e^{\otimes 2k} \right>
\nn \\
& = \lambda_k 
\left( e^{i_0}_{i_1} - \tfrac{1}{2} \, \delta^{i_0}_{i_1} \right) 
\otimes e^{i_1}_{i_2} \otimes e^{i_2}_{i_3} \cdots \otimes e^{i_{2k}}_{i_0} .
\end{align}
Here $\delta_{ij}$ is the usual Kronecker symbol and the $\lambda_k$ normalization constants.

Similarly, for a unitary 
$$
U = (U^{i}_{j}) \in \Mat_r (\ca) \, \qquad U U^* = U^* U,  
$$
the component $\ccc_{k+\frac{1}{2}} (U)$ of the Chern character of $U$ is the element of
$\ca \otimes (\wt{ \ca} )^{\otimes {(2k+1)}}$ given explicitly by the formula 
\begin{align}\label{B2}
\ccc_{k+\frac{1}{2}} (U) &= \left< \underbrace{U \otimes U^* \otimes U \otimes U^* \otimes \cdots U \otimes U^*}_{2(k+1)} 
- \underbrace{U^* \otimes U \otimes U^* \otimes U \cdots \otimes U^* \otimes U}_{2(k+1)} \right>  \nn \\
&= \lambda_k  
\left( U^{i_0}_{i_1} \otimes U^*{}^{i_1}_{i_2} \otimes U^{i_2}_{i_3} \otimes \dots \otimes U^*{}^{i_{2k+1}}_{i_0}
- U^*{}^{i_0}_{i_1} \otimes \dots \otimes U^{i_{2k+1}}_{i_0} \right) 
\end{align}
with $\lambda_k$ again normalization constants.

The crucial property of the
components $\ccc_k (e)$ or $\ccc_{k+\frac{1}{2}} (U)$ is that they define a {\it cycle} in the $(b,B)$
bicomplex of cyclic homology \cite{co85}, \cite{lo98}, that is 
\begin{equation}\label{B3}
B \ccc_k (e) = b \ccc_{k+1} (e) \qquad \mbox{or} \qquad 
B \ccc_{k+\frac{1}{2}} (U) = b \ccc_{k+\frac{3}{2}} (U) 
\end{equation}
where $b$ is the Hochschild boundary operator and $B$ is the Connes boundary operator.  

For a noncommutative spherical manifold \cite{cl01}, \cite{cdv02}, one asks that the components of the character vanish but a top one that, due to \eqref{B3} is then a (non zero) Hochschild cycle and plays the role of the volume form for the noncommutative manifold. Specifically, in even dimensions, for $n=2m$ one asks
\beq\label{B4}
\ccc_k (e) = 0 , \qquad \mbox{for all} \quad k = 0,1, \dots m -1,
\eeq
and $\ccc_m (e)$ (with $b \ccc_m (e) = 0$ from \eqref{B3}) is the volume form. Similarly, in odd dimensions, 
for $n=2m+1$ the vanishing condition becomes 
\beq\label{B5}
\ccc_{k+\frac{1}{2}} (U) = 0 , \qquad \mbox{for all} \quad k = 0,1, \dots m -1,
\eeq
and $\ccc_{m+\frac{1}{2}} (U)$ (with $b \ccc_{m+\frac{1}{2}} (U) = 0$ from \eqref{B3}) is the volume form.


\end{document}